\documentclass[3p,times]{article}
\usepackage{amsmath}
\usepackage{amsfonts}
\usepackage{mathrsfs}
\usepackage{amssymb}
\usepackage{times}
\usepackage{subcaption}
\usepackage{authblk}
\usepackage{float}
\usepackage{graphicx}%
\usepackage{chngcntr}
\newtheorem{theorem}{Theorem}
\newtheorem{proposition}[theorem]{Proposition}
\newtheorem{remark}[theorem]{Remark}
\newtheorem{exam}[theorem]{\bf Example}
\newtheorem{corollary}[theorem]{Corollary}
\newtheorem{lemma}[theorem]{Lemma}

\title{Perturbations of copulas and Mixing properties}
\author[1]{Martial Longla\thanks{mlongla@olemiss.edu}}
\author[2]{Fidel Djongreba Ndikwa\thanks{fidel.djong@yahoo.fr}}
\author[1]{Mathias Muia Nthiani.\thanks{mnmuia@go.olemiss.edu}}
\author[3]{Patrice Takam Soh.\thanks{ptakamsoh@yahoo.fr}}
\affil[1]{University of Mississippi, Department of mathematics}
\affil[2]{University of Maroua, Department of mathematics}
\affil[3]{University of Yaounde I, Department of mathematics}
\counterwithin{equation}{section}
\counterwithin{theorem}{subsection}
\begin{document}
\maketitle

\begin{abstract}
This paper explores the impact of perturbations of copulas on the dependence properties of the Markov chains they generate. We consider Markov chains generated by perturbations of copulas. Results are provided for the mixing coefficients $\beta_n$, $\psi_n$ and $\phi_n$. Several results  on mixing for the considered perturbations are provided. New copula functions are provided in connection with perturbations of variables that induce other types of perturbation of copulas not considered in the literature.
\end{abstract}



\section{Introduction}
In many areas of research, fitting an appropriate copula to real data is one of the major problems. Several authors have provided various families of copulas designed with specific properties. Longla (2014) proposed a set of copula families that generate exponential $\rho-$mixing Markov chains, to complete the list of other authors cited therein. It's been shown in Longla and Peligrad (2012) and by many other authors that most of the dependence and mixing coefficients are heavily influenced by copulas of the model of interest. The mixing coefficients of interest in this paper were explored by Bradley (2007). The $\psi$-mixing condition has its origin in the paper by Blum and al. (1963). They studied a different condition (“$\psi$*-mixing”) similar to this measure of dependence. It took its present form in the paper of Philipp (1969). For examples of mixing sequences, see Kesten and O'Brien (1976). The general definitions of these mixing coefficients are as follows. Given any $\sigma$-fields $\mathscr{A}$ and $\mathscr{B}$ and a defined probability measure $P$,
$$\beta(\mathscr{A},\mathscr{B})=\mathbb{E}\sup_{B\in \mathscr{B}}|P(B|\mathscr{A})-P(B)|,$$ 
$$\psi(\mathscr{A},\mathscr{B})=\sup_{B\in \mathscr{B}, A\in \mathscr{A}, P(A)\cdot P(B)>0 }\frac{P(A\cap B)-P(A)P(B)}{P(A)P(B)},$$
$$\phi(\mathscr{A},\mathscr{B})=\sup_{B\in \mathscr{B}, A\in \mathscr{A}, P(A)>0}|P(B|A)-P(B)|.$$ 
In applications, knowing approximately a copula $C(u,v)$ appropriate to the model of the observed data, minor perturbations of $C(u,v)$ are considered to look for a better fit for the data. Komornik and al. (2017) remark this and cite the class of perturbations introduced by Mesiar and al. (2015). We are investigating in this work the impact of such perturbations on $\beta$-mixing, $\psi$-mixing and $\phi$-mixing for stationary ergodic Markov chains in the setup of Longla (2013), Longla (2014) and Longla (2015).

Perturbations that we consider in this work have been studied by many authors. The closest work considering the impact of perturbations is that of Komornik and al. (2017) and was limited to dependence coefficients of two random variables such as Kendall's $\tau$, Spearman's $\rho$, Blomqvist's $\beta$ and tail dependence coefficients. Mesiar and al. (2015) looked at the perturbations of bivariate copulas via modification of the formula of the copula itself, while Sheikhi and al. (2020) looked at the perturbations of copulas via modification of the random variables that the copulas are used to represent the dependence structure of. Namely, they perturbed the copula of $(X,Y)$ by looking at the copula of $(X+Z, Y+Z)$ for some $Z$ independent of $(X,Y)$ that can be considered as noise. This case is also looked at in this work from a different perspective.
Mesiar and al. (2019) worked on perturbations induced by a modification of one of the random variables of the pair. Namely, the copula of $(X,Y)$ was perturbed to obtain the copula of $(X+Z, Y)$. In this work, we look at other perturbations of the vector $(X,Y)$ among other things.
One of the main results of this work is the link we establish between perturbations of copulas and mixing properties of the Markov chains that they generate. To the best of our knowledge, this setup has not been considered in the literature.  To do this, we shall define first each of the notions of interest. Before heading into copulas, perturbations and mixing, let us formulate the following useful result.
\begin{proposition}\label{sequence}
Suppose a sequence of complex numbers $a_n$ converges to $a$. For $0<\theta<1$,

If \quad $\displaystyle b_n=\sum_{i=1}^{n}{n \choose i}\theta^i(1-\theta)^{n-i}a_i$, then $\displaystyle \lim_{n\to\infty}b_n = a$.
\end{proposition}

The rest of this paper is structured as follows. In Subsection 1.1. we propose an overview of copulas and copula-based Markov chains including some important facts on copula densities. In Subsection 1.2 we propose an introduction to perturbations and their relationship with copulas. We emphasize on the fact that some of these perturbations can be presented as convex combinations of copulas, which makes it easier to analyze. Here, we propose a new result on joint distributions of $(X_0, X_n)$ for a perturbed Markov chain. A Lemma is proposed for commuting copulas; and copulas are derived for several types of perturbations including the Hoeffding upper bound and the independence copula. In Subsection 1.3 we present mixing coefficients and some known results relating them to copulas. In Subsection 1.4 we present the relationship between perturbations and the dependence coefficients. A proposition is provided for the copula densities  $m_1, m_2, m_3, m_4$ defined in Longla (2013). This subsection ends by a theorem on dependence coefficients and two examples of perturbations using the Independence copula and the Hoeffding upper bound. Section 2 is devoted to the relationship between perturbations and mixing properties of the Markov chains generated by copulas. Results are proposed for $\beta$-mixing, $\phi$-mixing and $\psi$-mixing. In Section 3 few graphs are proposed to justify some of the findings of this paper. Section 4 presents subsections devoted to proofs of the main results of this paper.

\subsection{Copula theory and Markov chains}
A 2-copula is a bivariate function $C:[0,1]\times[0,1]\to[0,1]$, such that $C(0,x)=C(x,0) = 0$, $C(1,x)=C(x,1) =x$ for all $x\in [0,1]$ and for all $[x_1,x_2]\times[y_1,y_2] \subset [0,1]^2,$ $C(x_1,y_1)+C(x_2,y_2)-C(x_1,y_2)-C(x_2,y_1)\ge 0.$ It's known that any convex combination of 2-copulas ($C(x,y)=a_1 C_1(x,y)+\cdots+a_k C_k(x,y),$ for any positive integer $k$ and any positive numbers $a_1, \cdots, a_k$ such that $ a_1+\cdots+a_k=1$) is a 2-copula (see Longla (2015)). Sklar's theorem states that if $X$ and $Y$ are random variables with joint distribution $H$ and marginal distributions $F$ and $G$ respectively, then the function $C$ defined by $C(F(x),G(y)) =H(x,y)$ for all $x,y\in \mathbb{R}$ is a 2-copula. This copula, called the copula of $(X,Y)$, captures the dependence structure of the two variables and is unique if the variables are continuous. In copula theory there are very famous classes of copulas with various properties. Among such copulas are the Hoeffding upper and lower bounds ($M(x,y)=\min(x,y)$ and $W(x,y)=\max(x+y-1, 0)$) and the independence copula $\Pi (x,y)=xy$. The Hoeffding upper bound characterizes perfect positive correlation and the Hoeffding lower bound characterizes perfect negative correlation, while the independence copula characterizes independence of the random variables. The later means that two random variables are independent if and only if their copula is $\Pi(u,v)$.

It's been shown in the literature (see Darsow and al. (1992) and the references therein) that if $(X_1, \cdots, X_n)$ is a Markov chain with consecutive copulas $(C_1, \cdots, C_{n-1})$, then the fold product given by 
$$C(x,y)=C_1*C_2 (x, y)=\int^1_0 C_{1,2}(x, t)C_{2,1}(t, y)dt$$ is the copula of $(X_1,X_3)$ and the $\star$-product given by
$$ C(x,y,z)=C_1\star C_2 (x, y,z)=\int_0^y C_{1,2}(x, t)C_{2,1}(t, z)dt$$ is the copula of $(X_1,X_2,X_3)$. The $n$-fold product of $C(x,y)$ denoted $C^n(x,y)$ is defined by the recurrence $C^{1}(x,y)=C(x,y)$, $$C^{n}(x,y)=C^{n-1}*C(x,y).$$ In this work, $A_{,i}$ denotes the derivative with respect to the variable at the position $i$. It is also clear and well-known that derivatives of copulas with respect to their variables exist almost everywhere. Also, if the Marginal distribution of the stationary Markov chain is the uniform distribution, then $C_{i,1}(x,y)$ is the transition probability from $X_i=x$ to the set $[0,y]$ (see Longla and Peligrad (2012) and the references therein). 
Among other things, in Longla (2013) was used but not stated a fact that can be generalized to the following result.
\begin{proposition}
Let $c: [0,1]\times [0,1]\to \mathbb{R}$ be a function such that $c(x,y)\ge 0$ for all $x,y\in [0,1]$. If $$\int_{0}^{y}\int_{0}^{1}c(t,s)dtds=y ,\quad \mbox{and}\quad \int_{0}^{x}\int_{0}^{1}c(s,t)dtds=x,$$ then $c(x,y)$ is the density of an absolutely continuous copula defined by
$$C(x,y)=\int_{0}^{x}\int_{0}^{y}c(s,t)dtds.$$
\end{proposition}
We can easily deduce from this proposition the following statement. 
\begin{corollary}
Let $c: [0,1]\times [0,1]\to \mathbb{R}$ be a function such that $c(x,y)\ge 0$ for all $x,y\in [0,1]$. If $$\int_{0}^{1}c(t,s)dt=1 ,\quad \mbox{and}\quad \int_{0}^{1}c(s,t)dt=1,$$ then $c(x,y)$ is the density of an absolutely continuous copula defined by
$$C(x,y)=\int_{0}^{x}\int_{0}^{y}c(s,t)dtds.$$
\end{corollary}
This corollary was applied  to conclude that the functions $m_1, \cdots, m_4$ are copula densities but not proven in Longla (2013).
\begin{exam}
Given any functions $h: [0,1] \rightarrow \mathbb{R},$ and $g: [0,1]\rightarrow \mathbb{R} $, if we set $a_1=\sup h$, $a_2=\sup g$, $b_1=\inf h$ and $b_2=\inf g$, then the following functions are copula densities.
\begin{enumerate}
\item 
$\displaystyle m_1(x,y)=\frac{a_2-g(x)h(y)+h(y)||g||_1 +g(x)||h||_1}{a_2+||g||_1||h||_1}$

\item $\displaystyle m_2(x,y)=\frac{a_1a_2-g(x)h(y)+h(y)||g||_1 +g(x)||h||_1}{a_1a_2+||g||_1||h||_1}$

\item $\displaystyle m_3(x,y)=\frac{a_2(a_1-b_1)-g(x)(a_1-h(y))+(a_1-h(y))||g||_1 +g(x)(a_1-||h||_1)}{a_2(a_1-b_1)+||g||_1 (a_1-||h||_1)}$

\item $\displaystyle m_4(x,y)=\frac{(a_2-b_2)(a_1-b_1)-(a_2-g(x))(a_1-h(y))}{(a_2-b_2)(a_1-b_1)+(a_2-||g||_1) (a_1-||h||_1)}+$ 

$\displaystyle +\frac{(a_1-h(y))(a_2-||g||_1) +(a_2-g(x))(a_1-||h||_1)}{(a_2-b_2)(a_1-b_1)+(a_2-||g||_1) (a_1-||h||_1)}$.
\end{enumerate}
\end{exam}
We also consider the Frank copula for permutations. This copula, investigated by De Baets and De Meyer (2013) (they proved that any copula of this family satisfies a set of inequalities good for stochastic interpretation in applications) is defined as follows:
$$C(u,v)=-\frac{1}{\lambda}\ln (1+\frac{(e^{-\lambda u}-1 )(e^{-\lambda u}-1)}{e^{-\lambda}-1}), \quad \lambda\ne 0.$$

Stationary Markov chains are defined by a marginal distribution $F$ and a copula $C$ that captures de strength of the dependence between any two consecutive variables of the chain. When the copula of the stationary Markov chain is symmetric, the chain is said to be reversible; this means that  the distribution of the time reversed chain is same as that of the original Markov chain and the distribution of any set of variables along the chain is invariant to time shifts.

\subsection{Perturbations of copulas}
Perturbations of copulas are tools used in applications to obtain better fits of real world data. They technically bring new families of copulas to the investigation. We consider the 2-copula $C$ and the related copula $C_H$ defined by $C_H(x,y)=C(x,y)+H(x,y),$
where $H : [0, 1]^2 \to \mathbb{R}$ is a continuous function. The function $H$ is the perturbation factor and the copula $C_H$ is called a perturbation of $C$. This perturbation method introduced by Durante and al. (2013) and its special case based on $$H(x,y)=H_\theta(x,y):=\theta (x-C(x,y))(y-C(x,y)), \quad\quad \theta\in [0,1]$$ was investigated by Mesiar and al. (2015).
A famous example of perturbation that follows this rule is the Farlie-Gumbel-Morgenstern copula family defined by $$C_\theta(x,y)=xy+\theta xy(1-x)(1-y), \quad \theta \in [0,1].$$ 
It is clear that $C_{\theta}(x,y)=\Pi (x,y)+\theta (x-\Pi(x,y))(y-\Pi(x,y))$. This shows that Farlie-Gumbel-Morgenstern copulas are perturbations of the independence copula. Moreover, it is clear that the perturbation with parameter $\theta$ is a convex combination of the original copula and the perturbation with parameter $1$. That is $$C_{H_\theta}(x,y)=C_\theta(x,y)=(1-\theta) C(x,y)+\theta C_{1}(x,y).$$ This fact is very valuable for the study of the mixing structure of Markov chains generated by the perturbation of $C$.
We also consider the perturbation of the form 

\begin{equation}
\hat{C}_\theta (u,v)=C(u,v)+\theta(M(u,v)-C(u,v)), \quad \theta\in [0,1], \label{CM}
\end{equation}

\begin{equation}
\tilde{C}_\theta (u,v)=C(u,v)+\theta(\Pi(u,v)-C(u,v)), \quad \theta\in [0,1], \label{CP}
\end{equation} 
and
\begin{equation}\label{Dolati}
C_1(u,v)=C(u,v)(u+v-C(u,v)).
\end{equation}
The copula defined by \eqref{Dolati} was introduced in Example 1.3.8 of Durante and Sempi (2015) as an example of copula generated from known copula families presented in Dolati and Ubeda (2009). To understand the stochastic interest of this copula, let $(\Omega , F , P )$ be a probability space, $( X_1 , Y_1 )$ and $(X_2 , Y_2)$ be two independent vectors from a population $(X,Y)$ with equal joint distribution functions given by the copula $C(u,v)$ . Let ( $Z_{(1)} ,Z_{(2)}$) denote the vector of the order statistics of any vector $(Z_1, Z_2)$. Consider the random vector $( Z_1 , Z_2 ) = ( X_{(1)} , Y_{(2)} )$ with probability $1 / 2$ and , $( X_{(2)} , Y_{(1)} )$ otherwise. The distribution function of $( Z_1 , Z_2 )$ is given by $C_1(u,v)$ defined by formula \eqref{Dolati}.

These perturbations are very different from $C_\theta$ obtained with $H_\theta$ given above. This is because $C_\theta$ is a non-linear function of the original copula $C(u,v)$ while \eqref{CM} and \eqref{CP} are linear functions of $C(u,v)$. This is very important when it comes to the impact of perturbations on dependence coefficients of the random variables with copula $C(u,v)$. The following holds for these perturbations.
\begin{proposition}\label{convex}
The joint distribution of $X_0$ and $X_n$, when $(X_1, \cdots, X_n)$ is a Markov chain generated by $\hat{C}_\theta(u,v)$ or $(\tilde{C}_\theta(u,v))$ is given respectively by 
\begin{equation}
\hat{C}_\theta^n(u,v)=\sum_{i=1}^{n}{n \choose i}\theta^{n-i} (1-\theta)^{i}C^{i}(u,v)+\theta^nM(u,v) \label{Mcop}
\end{equation} 
\begin{equation}
\mbox{or}\quad \tilde{C}_\theta^n(u,v)=(1-\theta)^nC^{n}(u,v)+(1-(1-\theta)^n)\Pi(u,v). \label{Pcop}
\end{equation} 
\end{proposition}
This proposition shows that along a Markov chain generated by a perturbation copula, the copulas remain in the class of given perturbations when the independence copula is used.
Proposition \ref{convex} relies on the following result.
\begin{lemma} \label{commut}
For any  $2$-copulas $A(u,v)$ and $B(u,v)$ that commute, the copula 

$C(u,v)=\theta A(u,v)+(1-\theta)B(u,v)$ satisfies the following equation
\begin{equation}
C^{n}(u,v):=(\theta A+(1-\theta)B)^n(u,v)=\sum_{k=0}^{n}{n \choose k}\theta^k(1-\theta)^{n-k}A^k*B^{n-k}(u,v).
\end{equation}
\end{lemma}

Masiar and al. (2019) have discussed the form of perturbation of a copula $C_{X+Z,Y}$ when the first random variable $X$ is affected by a noise $Z$ (independent of both $X$ and $Y$), leaving as open question the case of the perturbation of both variables $X$ and $Y$ by $Z_1$ and $Z_2$. Dropping the independence of $Z$ or ($Z_1, Z_2$) and $X$ or $(X,Y)$ and generalizing to $n$–dimensional random vectors affected by noise were left for future research. We consider in this paper the case of Markov chains.

\begin{proposition} \label{independent}
Let $X,Y,Z$ be such that $(X,Y)$ has copula $C(u,v)$ and joint distribution $H$.
Assume that $X, Y, Z$ have respectively cumulative distributions $F_1$, $F_2$ and $F_3$. 
Moreover, assume that $Z$ is independent of $(X,Y)$. Then, the copula of $(X+Z, Y)$ is found by the formula
\begin{equation}\label{newcopula}
C_5 (u,v)=\int_0^1 C(F_1(F^{-1}_4(u)-F_3^{-1}(t)), v)dt, 
\end{equation}
where $\displaystyle F_4(x)=\int_{-\infty}^{\infty}F_1(x-t)f_3(t)dt$ is the cumulative distribution function of $W=X+Z$ and $F^{-1}(t)=\inf\{x: F(x)=t\}$.
\end{proposition}
This proposition shows that any choice of three cumulative distributions will give us a new copula function via formula \eqref{newcopula}. The fact that formula \eqref{newcopula} defines a copula follows from Sklar's theorem (See Sklar (1959) for more on the topic). The other assertion of the proposition will be proven in the section devoted to proofs. On the other hand, if one wants to show that the function is grounded (meaning $C_5(u,0)=C_5(0,v)=0$), this is easily obtained in this case. $C(F_1(F^{-1}_4(u)-F_3^{-1}(t)), 0)=0$ because $C$ is a copula and for any value of $t$, and any value of $v$, we have $$C(F_1(F^{-1}_4(0)-F_3^{-1}(t)), v)=C(F_1(-\infty), v)=0,$$ using the fact that $F^{-1}_4(0)=-\infty$ and $C$ is a copula. It is also easy to show that random variables with cumulative distribution $C_5$ have uniform marginals, because for all $v,t$ 
$$C(F_1(F^{-1}_4(1)-F_3^{-1}(t)), v)=C(F_1(\infty), v)=C(1,v)=v.$$
The last equality of the boundary conditions is $C_5(u,1)=u$. This follows from the fact that
$$C_5 (u,1)=\int_0^1 C(F_1(F^{-1}_4(u)-F_3^{-1}(t)), 1)dt=\int_0^1 F_1(F^{-1}_4(u)-F_3^{-1}(t))dt,
$$
from which we have $$C_5(u,1)=\int_{-\infty}^{\infty}F_1(F^{-1}_4(u)-t)f_3(t)dt=F_4(F^{-1}_4(u))=u.$$ 
Notice that such a perturbation doesn't impact the independence copula. This is obviously true because if we begin with a set of three independent random variables $X,Y,Z$, then $X+Z$ and $Y$ are independent and their copula is the independence copula. One of the questions is to study the impact of this transformation of the copula $C$ on the mixing structure of the Markov chains it generates. 
\begin{exam}[Examples of Perturbations]

Considering some famous copulas, we have the following examples of copulas, obtained from perturbations of the type \eqref{newcopula}.
\begin{enumerate}
\item If $(X,Y,Z)$ are independent, then $C(u,v)=uv$ and the copula of $(X+Z,Y)$ is $C_5(u,v)=uv$.
\item If $(X,Y)$ has distribution $M(u,v)$ and $Z$ is Uniform on $(0,1)$. Then the copula of $(X+Z, Y)$ is

$C_5(u,v)=
\begin{cases}
u, & \text{if\quad $\sqrt{2u} \leq v\leq 1 $} \\
v\sqrt{2u}-\frac{1}{2}v^2,& \text{if\quad $v \leq \sqrt{2u} \leq 1 $}\\
v-\frac{1}{2}(1-v-\sqrt{2(1-u)})^2,& \text{if\quad $ 0 \leq 1-\sqrt{2(1-u)}\leq v$}\\
v, & \text{if\quad $v\leq 1-\sqrt{2(1-u)}$.}
\end{cases}
$

A Markov chain generated by this copula started at any point other than 0 or 1 can reach any subset in one or more steps.\\
- When started at $u>.5$, it can reach any point above $1-\sqrt{2(1-u)}$ (this number is below $u$). \\
- Thus, the set of states that can be reached is always enlarging, getting closer to $.5$, until the first number below $.5$ is reached. And from any number below $.5$, the chain will easily reach any number.\\
The lower tail dependence coefficient for the copula $M(u,v)$ is $\lambda_L(M)=1$, but the lower tail dependence coefficient for $C_5(u,v)$ is $\lambda_L(C_5)=0$. In this sense, this method of perturbation is stronger than $C_H$. It turns a tail dependent copula into one that exhibits no tail dependence. More copulas can be generated in this case using different Marginal distributions for $Z$.

The upper tail dependence coefficient is $\lambda_U(C_5)=0$, while $\lambda_U(M)=1$.

\end{enumerate}
\end{exam}

\begin{proposition}\label{independentboth}
Let $X,Y,Z$ be such that $(X,Y)$ has copula $C(u,v)$ and joint distribution $H$.
Assume that $X, Y, Z$ have respectively cumulative distributions $F_1$, $F_2$ and $F_3$. 
Moreover, assume that $Z$ is independent of $(X,Y)$. Then, the copula of $(X+Z, Y+Z)$ is found by the formula
\begin{equation}\label{newcopula}
C_6 (u,v)=\int_0^1 C(F_1(F^{-1}_4(u)-F_3^{-1}(t)), F_2(F^{-1}_5 (v)-F_3^{-1}(t)))dt, 
\end{equation}
where $\displaystyle F_5(x)=\int_{-\infty}^{\infty}F_2(x-t)f_3(t)dt$ is the cumulative distribution function of $W=Y+Z$, $F_4$ is defined in Proposition \ref{independent} and $F^{-1}(t)=\inf\{x: F(x)=t\}$.
\end{proposition}

Proposition \ref{independentboth} exhibits a perturbation that modifies the independence copula. An alternative form of the copula $C_6(u,v)$ was given in Sheikhi and al. (2020). Here we provide a simple integral based on computing probabilities via conditioning while they had provided a double integral to find the given copula. 
\begin{exam}[Double perturbation] \label{doublepert}

Considering the independence case and the extreme case of $M(u,v)$, we have the following:
\begin{enumerate}
\item If $(X,Y,Z)$ are independent, then $C(u,v)=uv$ and the copula of $(X+Z,Y+Z)$ is 
\begin{equation}
C_6(u,v)=\int_0^1 F_1(F^{-1}_4(u)-F_3^{-1}(t))F_2(F^{-1}_5 (v)-F_3^{-1}(t))dt.
\end{equation}
It is clear in this case, that unlike for the perturbation in Proposition \ref{independent}, there are many copulas coming out of this representation depending on the marginal distributions of the original random variables. If we consider uniform distributions for all three variables, then the copula is

$C_6(u,v)=
\begin{cases}
\frac{1}{2}u\sqrt{2v}, & \text{if\quad A} \\
u-\frac{1}{6}(\sqrt{2u}+\sqrt{2(1-v)}-1)^3,& \text{if\quad B}\\
u,& \text{if\quad C}\\
u-\frac{1}{3}(1-v)(\sqrt{2(1-v)}+3-3\sqrt{2(1-u)}) & \text{if\quad D}\\
v-\frac{1}{3}(1-u)(\sqrt{2(1-u)}+3-3\sqrt{2(1-u)}) & \text{if\quad E}\\

v-\frac{1}{6}(\sqrt{2v}+\sqrt{2(1-u)}-1)^3,& \text{if\quad F}\\
v,& \text{if\quad G}\\

\frac{1}{2}v\sqrt{2u}, & \text{if\quad H}, \\
\end{cases}
$

where $A=\{u\leq v \leq 1/2\} $, $B=\{1/2< v\leq 1-\frac{1}{2}(\sqrt{2u}-1)^2, u\leq 1/2\}$, 

$C=\{v> 1-\frac{1}{2}(\sqrt{2u}-1)^2, u\leq 1/2\}$,  $D=\{1/2\leq u\leq v \}$, 

$E=\{1/2\leq v\leq u \}$,  $F=\{1/2< u\leq 1-\frac{1}{2}(\sqrt{2v}-1)^2, v\leq 1/2\}$, 

$G=\{u> 1-\frac{1}{2}(\sqrt{2v}-1)^2, v\leq 1/2\}$ and $H=\{v\leq u \leq 1/2 \}$.

\item If $(X,Y)$ has distribution $M(u,v)$ and $Z$ has any distribution, then the copula of $(X+Z, Y+Z)$ is $M(u,v)$.
In this case, we see that the perturbation doesn't affect the Hoeffding upper bound. But with different Marginal distributions for $X$ and $Y$, keeping their copula as $M(u,v)$, the effect is obvious.

A Markov chain generated by this copula is aperiodic and irreducible (thus, Harris recurrent).
The lower tail dependence coefficient for the copula $M(u,v)$ is $\lambda_L(\Pi)=0$ and the lower tail dependence fcoefficient for $C_6(u,v)$ is $\lambda_L(C_6)=0$. The upper tail dependence coefficient is $\lambda_U(C_6)=0$ and $\lambda_U(\Pi)=0$. Thus, for this method, perturbations of $\Pi$ don't exhibit tail dependence.

\end{enumerate}
\end{exam}

Consider a perturbation in which both variables are modified by two different independent random variables. 
\begin{proposition} \label{bothcomp} Assume $(X,Y)$ is independent from $(Z_1, Z_2)$. Also, assume $Z_1, Z_2$ are independent and the copula of $(X,Y)$ is $C(u,v)$. If variables $X,Y, Z_1$, $ Z_2$, $X+Z_1$ and $Y+Z_2$ have cumulative distributions $F_1, F_2, G_1$, $G_2$, $F_7$ and $F_8$ respectively, then the copula of $(X+Z_1, Y+Z_2)$ is
\begin{equation}
C_7(u,v)=\int_{-\infty}^{\infty}\int_{-\infty}^{\infty}C(F_1(F_7^{-1}(u)-t_1), F_2(F_8^{-1}(v)-t_2))dG_2(t_2)dG_1(t_1). \label{copulamixt}
\end{equation}

When all variables are uniform on $(0,1)$, this formula gives
\begin{equation} C_7(u,v)=
\begin{cases}
\int_0^{b(u)}\int_0^{b(v)}C(b(u)-t_1, b(v)-t_2)dt_2dt_1, & \text{if\quad I} \\
\int_0^{b(u)}\int_{a(v)}^{1}C(b(u)-t_1, 1+a(v)-t_2)dt_2dt_1+ua(v), & \text{if\quad II}\\
\int_{a(u)}^{1}\int_0^{b(v)}C(1+a(u)-t_1, b(v)-t_2)dt_2dt_1+va(u), & \text{if\quad III}\\
\int_{a(u)}^{1}\int_{a(v)}^1 C(1+a(u)-t_1, 1+a(v)-t_2)dt_2dt_1+\\+ua(v)+va(u)-a(u)a(v), & \text{if\quad IV},
\end{cases}
\end{equation}
where 
$$a(x)=1-\sqrt{2(1-x)}, \quad b(x)=\sqrt{2x},$$ $$I=\{u\leq 1/2, v \leq 1/2 \},\quad II=\{ u\leq 1/2, v \ge 1/2\},$$ 
$$III=\{v\leq 1/2, u\ge 1/2\}\quad and \quad IV=\{v\ge 1/2, u\ge 1/2\}.$$
\end{proposition}

Notice that formula \eqref{copulamixt} can be rewritten as 
\begin{equation}
C(u,v)=\int_{0}^{1}\int_{0}^{1}C(F_1(F_7^{-1}(u)-G_1^{-1}(u_1)), F_2(F_8^{-1}(v)-G_2^{-1}(u_2)))du_2du_1. \label{copulamixt1}
\end{equation}

\begin{exam}[Double perturbation with different independent variables]
Let’s consider a few cases of the setup of the copula defined by \eqref{copulamixt}.
\begin{enumerate}

\item If $(X,Y)$ are independent, then for any distributions of the random variables, $C_7(u,v)=uv$. This is true because under this scenario, $(X,Z_1)$ and $(Y,Z_2)$ are independent. This can also be easily checked using the given formulas to compute the integrals. It is an easy exercise.
\item If $(X,Y)$ have copula $M(u,v)$, then the copula $C_7(u,v)$ is same as that of a perturbation of independent random variables by a common variable $Z$. This is true because the copula M describes a relationship $Y=Z$ almost everywhere. It can also be checked using directly the formulas.
\end{enumerate}

\end{exam}
\subsection{Mixing coefficients for copula-based Markov chains}

Given the alternative form of the transition probabilities for a Markov chain generated by an absolutely continuous copula and a marginal distribution with positive density, it was shown by Longla and Peligrad (2012) that for $\mathscr{A}=\sigma(X_{i}, i\leq 0)$, $\mathscr{B}=\sigma(X_{i}, i\geq n)$ and the Lebesgue measure $\lambda$, these coefficients have the following form:
$$ \beta_{n}=\int_{0}^{1}\sup_{B\subset I}|\int_{B}(c_{n}(x,y)-1)dy|dx,$$
$$\psi_n=\sup_{B\in \mathscr{B}, A\in \mathscr{A}, \lambda(A)\cdot \lambda(B)>0 }\frac{\int_A\int_B c_n(x,y)-1dxdy}{\lambda(A)\lambda(B)},$$
$$\phi_{n}=\sup_{B\subset I}ess\sup_{x}|\int_{B}(c_{n}(x,y)-1)dy|,$$
where, $c_n$ is the density of the random vector $(X_{0}, X_{n})$. 

A stochastic process is said to be $\psi$-mixing, if $\psi_{n}\rightarrow 0$; $\beta$-mixing if $\beta_{n}\rightarrow 0$ and $\phi$-mixing if $\phi_{n}\rightarrow 0$. The process is exponentially mixing, if the convergence rate is exponential. A stochastic process is said to be geometrically ergodic, if $\beta_{n}$ converges to $0$ exponentially fast. For more about the relationship between copulas and stochastic processes see Darsow and al. (1992) and for results on mixing for copula-based Markov chains see Beare (2010), Longla and Peligrad (2012), Longla (2013), Longla (2014) and the references therein. For convex combinations of copulas and their mixing structure, see Longla (2015). These mixing conditions are very important in establishing central limit theorems for functions of the Markov chains generated by copulas in applications. Various limit theorems have been extended from sequences of independent random variables to random sequences satisfying one of the mixing conditions with some specified mixing rate; a mixing rate being any relationship that gives the rate of convergence of the mixing coefficient. See Phillip (1971) for limit theorems and Bradley (2005) for a survey on mixing sequences and for more on their applications, see Bradley (2007).

\subsection{Perturbations of copulas and dependence}
Several dependence coefficients are used in the literature for different reasons. Though we will focus more on mixing conditions in this paper, we find it useful to introduce the dependence coefficients to clarify the differences and known results. The most popular dependence coefficient of two random variables $X$ and $Y$ is the Pearson linear correlation coefficient $\rho(X,Y)$. Among others, there is the Kendall's $\tau$, the Spearman's $\rho_S$, the Blomqvist $\beta$, the upper tail dependence $\lambda_U$, the lower tail dependence $\lambda_L$ defined as follows in terms of copulas.
The Spearman's $\rho_S$ is defined as the difference between the probability of concordance and the probability of discordance of the two pairs of random variables $(X_1, Y_1)$ and $(X_2, Y_3)$ formed from 3 independent copies of the vector $(X,Y)$. If $C(u,v)$ is the copula of $(X,Y)$, then $$\rho_S(C) = 12 \int_{0}^1\int_{0}^1 C(u,v)dudv-3.\quad \text{Kendall's $\tau$ is given by} $$  $$\tau(C)=4\int_{0}^1\int_{0}^1 C(u,v)dC(u,v)-1=1-4\int_{0}^1\int_{0}^1 C_{,1}(u,v)C_{,2}(u,v)dudv.$$ Kendall’s $\tau$ is the difference between the probability of concordance and discordance of two copies of the vector $(X,Y)$. Blomqvist's $\beta$ in terms of copulas is given by $$\beta(C)=4C(1/2, 1/2)-1$$ and is called the medial correlation coefficient.
The Gini's $\gamma$ is the measure of concordance difference between $C(u,v)$ and copulas $M(u,v)$ and $W(u,v)$. It is given in terms of copulas by the formula $\displaystyle \gamma(C)=4\int_{0}^{1}C(u,u)+C(u,1-u)du.$
The upper tail dependence coefficient is given by $\displaystyle \lambda_U(C)=\lim_{u\to 1}\frac{1-2u+C(u, u)}{1-u}$ and the lower tail dependence is defined as $\displaystyle \lambda_L(C)=\lim_{u\to 0}\frac{C(u, u)}{u},$ (see Nelsen (2006) for more on these coefficients).
These coefficients play an important role in establishing the limit theorems when there is dependence among observations. Tail dependence coefficients are also used in risk management. 
Copulas with densities $m_1, \cdots, m_4$ generate Markov chains that exhibit $\rho$-mixing, but the main drawback is in the following proposition.
\begin{proposition}\label{mcopula}
Copulas with densities $m_1, \cdots, m_4$ do not exhibit lower nor upper tail dependence, no matter what are the functions used in their definition. 
\end{proposition}
For these copulas, we are interested in knowing what kind of dependence can be when they are perturbed. The influence of perturbations that we investigate has been studied by Komornik and al. (2017). They have shown that the difference between the values of the considered coefficients for the perturbation and the copula is proportional to the perturbation parameter $\alpha\in [0, 1]$ for the Spearman’s $\rho$, Blomqvist’s $\beta$ and Gini’s $\gamma$, while it is a linear combination of the form $c_1\alpha + c_2\alpha^2$ (with clearly dominative linear coefficient) for Kendall’s $\tau$. Moreover, they have shown that this class of perturbations does not change the values of the coefficients of tail dependence along the main diagonal, but yields their reduction along the second diagonal. Based on formulas of these coefficients, the following  can be obtained via simple computations.

\begin{theorem} For any random variables with copulas defined by \eqref{CP} or \eqref{CM}, the following holds:
\begin{enumerate}
\item $\rho_S(\tilde{C}_\theta)=(1-\theta)\rho_S(C)$ and $\rho_S(\hat{C})=(1-\theta)\rho_{S}(C)+\theta$;
\item $\gamma(\tilde{C}_\theta)=(1-\theta)\gamma(C)+2\theta$ and $\gamma(\hat{C})=(1-\theta)\gamma(C)+3\theta$;
\item $\beta(\tilde{C}_\theta)=(1-\theta)\beta(C)$ and $\beta(\hat{C})=(1-\theta)\beta(C)+\theta$;
\item $\lambda_U(\tilde{C}_\theta)=(1-\theta)\lambda_U(C)$ and $\lambda_U(\hat{C})=(1-\theta)\lambda_U(C)+\theta$;
\item $\lambda_L(\tilde{C}_\theta)=(1-\theta)\lambda_L(C)$ and $\lambda_L(\hat{C})=(1-\theta)\lambda_L(C)+\theta$.
\end{enumerate}
\end{theorem}

\section{Perturbations of copulas and mixing}
The connection between perturbations of copulas and mixing is first found in Longla (2015), where it is shown that the Farlie-Gumbel-Morgenstern family of copulas generates $\psi'$-mixing Markov chains. The mixing coefficient $\psi'$ is not defined here, but is mentioned because $\psi$-mixing implies $\psi'$-mixing and $\psi'$-mixing implies $\phi$-mixing. Thus, Markov chains generated by the Farlie-Gumbel Morgenstern copula family are $\phi$-mixing. Though it was not stated in Longla (2015) that these copulas are perturbations of the independence copula, it is clear that the mixing structure of the Markov chain generated by a perturbation copula is complex even in the case of the obvious lack of dependence in the original copula as shown by the theorems therein. It is our goal to establish the connection between perturbations, copulas and mixing in general for the considered mixing coefficients.
\subsection{On $\beta$-mixing}
$\beta$-mixing conditions have long been studied by many authors, including Beare (2010), Longla and Peligrad (2012), Longla (2014) and Longla (2015). Many other references can be found in the cited works. We do not provide here a complete list because this work is mostly about perturbations. It is well known that a strictly stationary Markov chain is absolutely regular if and only if it is irreducible and aperiodic (see Bradley (2005) or Longla and Peligrad (2012)). Longla (2015) provided results on mixing for convex combinations of copulas that we rely on in this work. Recalling that for copula-based Markov chains generated by an absolutely continuous copulas, $$ \beta_{n}:=\beta(C^n)=\int_{0}^{1}\sup_{B\subset I}|\int_{B}(c_{n}(x,y)-1)dy|dx,$$
if we denote $\beta(C)$ the coefficient of $\beta$-mixing for two random variables with copula $C$, we can easily prove the following: 
\begin{remark} In general, for any convex combination of two copulas (here $0\leq a\leq 1$), the $\beta$-mixing coefficient satisfies the following inequalities:
\begin{equation}
\beta(a C_1+(1-a)C_2)\leq a\beta(C_1)+(1-a)\beta(C_2), \label{beta1}
\end{equation}
\begin{equation}
\beta(a C_1+(1-a)C_2)\ge a\beta(C_1)-(1-a)\beta(C_2). \label{beta2}
\end{equation}
When $C_2(u,v)=\Pi(u,v)$, it follows from \eqref{beta1} and \eqref{beta2} that 
\begin{equation}
\beta(a C_1+(1-a)\Pi)= a\beta(C_1). \label{fbeta}
\end{equation}
\end{remark}
\begin{theorem}[$\beta$-mixing and perturbations] \label{theobeta} {\quad} {\quad}\\
\begin{enumerate}
\item Any perturbation of the form $\tilde{C}_{\theta}$ generates geometrically ergodic Markov chains.
\item Any perturbation of the form $\hat{C}_{\theta}$ generates absolutely regular Markov chains, if the original copula generates absolutely regular Markov chains.
\item Any perturbation of copulas of the form $C_{H_\theta}$ using $C_1(u,v)$ generates absolutely regular Markov chains when the original copula generates absolutely regular Markov chains.
\end{enumerate}
\end{theorem}

Thus, a perturbation of a copula by means of the independence copula rescales the $\beta$-mixing coefficient and increases exponentially the rate of convergence. In fact, it turns non-mixing sequences into geometrically ergodic sequences. Alternatively, we can say that a perturbation of the independence copula by means of any copula $C$ gives a copula that generates geometrically ergodic Markov chains. So, the selection of the perturbation copula is very important in general for the $\beta$-mixing property, but irrelevant when the independence copula is perturbed. The result for the case when the perturbation uses $M(u,v)$ is as follows. 
\begin{proposition}[Mixing for perturbations using $M(u,v)$] \label{Mmixing} {}\quad{}\\
For any copula $C(u,v)$ and ${\theta}\in(0,1)$, 
$\beta_n(\hat{C}_{\theta})\to c$, if $\beta_n(C)\to c$ as $n\to \infty$.
\end{proposition}

The above mentioned cases show that the dependence structure is seriously affected, but the mixing structure depends on the perturbation. In this sense it is beneficial to consider perturbations because they enlarge the class of central limit theorems that hold for the data.
Theorem 4.1 of Bradley (2005) states the following

\begin{theorem} \label{beta}
For a strictly stationary mixing Markov chain, either \begin{equation}\lim_{n\to\infty}\beta_n(C)=0 \quad or\quad \beta_n(C)=1,\quad for\quad all \quad n\ge 1.\end{equation}
\end{theorem}
\begin{remark}\label{rbeta}
Theorem \ref{beta} implies that to show that a Markov chain generated by the copula $C(u,v)$ is $\beta$-mixing, it is enough to show that $\beta(C)<1.$ This is true because $0\leq \beta(C)\leq 1$ by definition and the theorem states that failing to have $\beta$-mixing (or absolute regularity) is equivalent to having a constant sequence $\beta_n(C)=1$.
\end{remark}
Remark \ref{rbeta} and formula \eqref{fbeta} are enough to show that any convex combination of copulas containing the independence copula is geometric $\beta$-mixing (in other words geometrically ergodic). This proves that the perturbation of any copula by means of the independence copula creates a copula that generates geometrically ergodic Markov chains as stated above.

\subsection{On $\phi$-mixing}
A reformulation of Theorem 3.3 and Theorem 3.4 of Bradley (2005) for $\phi$-mixing is as follows.
\begin{theorem}\label{phi}
Suppose $(X_n,n\in Z)$ is a (not necessarily stationary) Markov chain. Then each of the following statements holds:
\begin{enumerate}
\item If $\phi_n<1/2$ for some $n\ge 1$, then $\phi_n\to 0$ at least exponentially fast as $n\to\infty$.
\item Suppose $(X_n,n\in Z)$ is ergodic and aperiodic. If $\phi_n<1$ for some $n\ge 1$, then $\phi_n\to 0$ (at least exponentially fast) as $n\to\infty$.
\end{enumerate}
\end{theorem}
Given that here we are not talking about densities, we cannot use the results of Longla (2015) on copulas with density bounded away from $0$.
Copulas in this work have a rather general form. We can formulate the following.
\begin{theorem}\label{phi}
Let $C(u,v)$ be a copula. Consider the perturbations $\tilde{C}_{\theta}(u,v)$ and $\hat{C}_{\theta}(u,v)$. The following holds.
\begin{enumerate}
\item if $\phi_n(C)\to c$ as $n\to\infty$, then $\phi_n(\hat{C}_{\theta})\to c$ as $n\to \infty$. 
\item $\phi_n(\tilde{C}_{\theta})\to 0$ exponentially fast as $n\to 0$.
\end{enumerate}
\end{theorem}

The first part of Theorem \ref{phi} shows that the perturbation will generate $\phi$ mixing Markov chains if the original copula generates $\phi$-mixing Markov chains. If we look at the formula as a perturbation of the Hoeffding upper bound copula by means of $C(u,v)$, this shows that the perturbation of $M(u,v)$ carries the $\phi$-mixing property of the copula that is used to perturb.  The second part of the theorem shows that any perturbation of the independence copula generates exponential $\phi$-mixing Markov chains.

A consequence of Theorem \ref{phi} is that any perturbation of the Frank copula using $M(u,v)$ or $\Pi(u,v)$ generates $\phi$-mixing Markov chains. This is true because the density of the Frank copula is bounded away from zero on $(0,1)^2$ (for instance, for $\lambda>0$, $c(u,v)>\lambda e^{-2\lambda}/(1-e^{-\lambda})$) and by Longla (2015) it generates $\phi$-mixing. 
\subsection{On $\psi$-mixing}
Kesten and O'Brien (1976) use the notation $\rho$ to denote the formula that we present here for the coefficient of $\psi$-mixing and use the name $\psi$-mixing sequence for sequences that we name so in this work. They provided a set of examples of mixing sequences that one can use. Bradley (1983) proposed a study of $\psi$-mixing sequences for stationary sequences and some references on the topic. The $\psi$-mixing coefficient he used differs from the one used here by an added term $1$. The coefficient used here measures the rate of convergence of the one used in Bradley (1983), and its convergence to $0$ implies the same consequences for the sequences of interest. It is good to mention that having a Markov chain makes the work easier, thanks to the use of the Markov property. A reformulation of his Theorem 1 is as follows
\begin{theorem}\label{psi}
For a strictly stationary mixing sequence, either $\psi_n=\infty $ for all $n$ or $\psi_n\to 0$ as $n\to \infty$. 
\end{theorem}
Based on this finding, if we want to show that a stationary Markov chain is $\psi$-mixing, it is enough to show that it is mixing and $\psi_1\ne \infty$.
A similar remark as that for $\beta$-mixing Markov chains holds.
\begin{remark}\label{psirem} In general, for any convex convolution of two copulas (here $0\leq a\leq 1$), the $\psi$-mixing coefficient satisfies the following inequalities:
\begin{equation}
\psi(a C_1+(1-a)C_2)\leq a\psi(C_1)+(1-a)\psi(C_2), \label{psi1}
\end{equation}
\begin{equation}
\psi(a C_1+(1-a)C_2)\ge a\psi(C_1)-(1-a)\psi(C_2). \label{psi2}
\end{equation}
When $C_2(u,v)=\Pi(u,v)$, it follows from \eqref{psi1} and \eqref{psi2} that 
\begin{equation}
\psi(a C_1+(1-a)\Pi)= a\psi(C_1).
\end{equation}
\end{remark}
We can now formulate for $\psi$-mixing a result for convex combinations of copulas in the spirit of the results of Longla (2015) for $\rho$-mixing sequences. 
\begin{theorem} \label{psimixing}
A convex combination of copulas generates stationary $\psi$-mixing Markov chains if each of the copulas of the combination generates $\psi$-mixing stationary Markov chains. 
\end{theorem}
Theorem \ref{psimixing} shows that perturbations by means of a copula $C_1(u,v)$ as described above do not influence the $\psi$-mixing property the way it does for dependence coefficients. Moreover, in the case of a perturbation by means of the independence copula, the convergence rate exponentially increases, while perturbations by means of the copula $M(u,v)$ are never $\psi$-mixing. This last statement is true because $\psi(M)=\infty$ and by Bradley's result $\psi_n((1-\theta) C+ \theta M)$ does not converge to $0$ for $\theta\in(0,1]$. Based on the perturbation of the form $C_\theta(u,v)$ and Theorem \ref{psimixing} we can state the following.
\begin{corollary}
The perturbation of a copula ($C_\theta(u,v)$) generates $\psi$-mixing stationary Markov chains if both the copula ($C(u,v)$) and the perturbation with parameter $1$ ($C_1(u,v)$) generate $\psi$-mixing stationary Markov chains.
\end{corollary}

\section{Graphs and final comments}

\begin{figure}[h]
\centering
\begin{subfigure}[b]{0.4\linewidth}
\includegraphics[width=\linewidth]{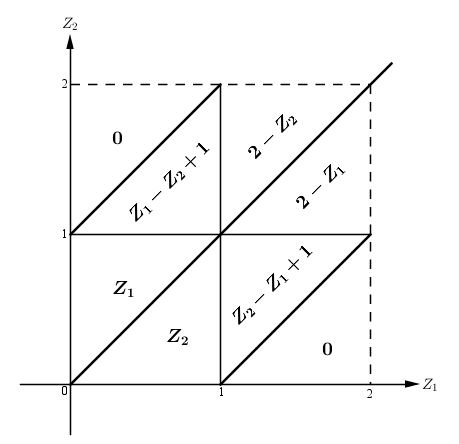}
\caption{Graph for joint probabilities} \label{Fig1a}
\end{subfigure}
\begin{subfigure}[b]{0.5\linewidth}
\includegraphics[width=\linewidth]{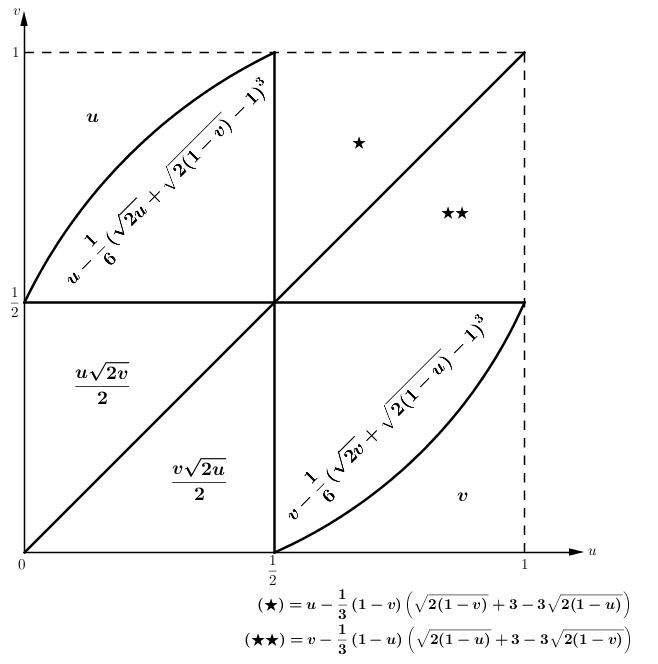}
\caption{Values of the copula of $(X+Z_1, Y+Z_2)$.} \label{Fig1b}
\end{subfigure}
\caption{Copula of $(X+Z_1, Y+Z_2).$}
\end{figure}

Figure \ref{Fig1a} presents the table of values of the conditional probability to use in each of the cases to find $P(X+Z_1<x, Y+Z_2<y)$, using the fact that when the variables of integration fall in each of the regions, the written function must be integrated. Figure \ref{Fig1b} presents the values of the copula of $(X+Z_1, Y+Z_2)$ when $(u,v)$ is in the indicated regions. This graph gives the copula of the perturbation in the case when $Z_1$ and  $Z_2$ are independent Uniform$(0,1)$ variables and $(X,Y)$ being independent from $(Z_1, Z_2)$ follows the distribution $M(u,v)$. It also represents the copula of the perturbation when $X$ and $Y$ are independent Uniform$(0,1)$ variables and $Z_1=Z_2=Z$ is uniformly distributed on $(0,1)$ and independent of $(X,Y)$. In the regions where the value of the copula is $u$ or $v$, the density of the copula is $0$. This means that points in this region can not be reached in one step for a Markov chain generated by this copula. It is clear from Figure \ref{Fig1b} that any point of the unit square can be reached in two steps. 

\begin{figure}[h!]
\centering 
\includegraphics[width=0.5\linewidth]{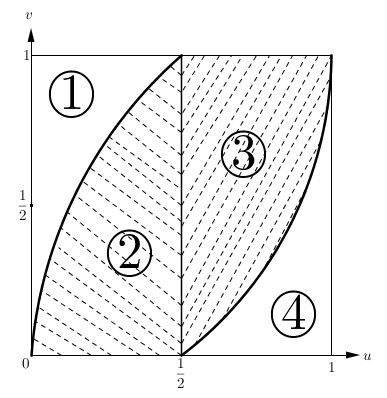}
\caption{Copula of $(X+Z, Y)$}
\label{fig2}
\end{figure}

Figure \ref{fig2} represents the regions used for computations of the copula of $(X+Z, Y)$ given above as $C_5(u,v)$. The shaded region is the range where the density of the copula exists and is not equal to $0$. In regions 1 and 4 the copula has density equal to $0$. Same as above, points in these two regions can not be reached in one step. It is good to mention that the copula $C_5(u,v)$ is not symmetric while the copula $C_6(u,v)$ is symmetric. Note that symmetric copulas generate reversible Markov chains.

\section{Proofs}
\subsection{Proof of Proposition \ref{convex}}
The proof of proposition \ref{convex} is an application of Lemma \ref{commut}. This is true because the copulas $M(u,v)$ and $\Pi(u,v)$ commute with any copula $C(u,v)$. Moreover, $$M*C(u,v)=C*M(u,v)=C(u,v)$$ and $$\Pi*C(u,v)=C*\Pi(u,v)=\Pi(u,v).$$ Formula \eqref{Mcop} is a direct application of formula \eqref{Pcop} and uses the fact that $$(1-\theta +\theta)^n=\sum_{i=0}^{n-1}{n \choose i}\theta^{n-i}(1-\theta)^i+(1-\theta)^n=1.$$
Lemma \ref{commut} follows from distributivity of the fold-product on convex combinations of copulas. The $n$-fold product here is a sum of $2^n$ terms that can be grouped by the number of times $A$ appears in them. Doing this, using commutativity, it is clear that all terms containing $k$ times $A$ are equal. Therefore, the binomial formula holds as 
$$(\theta A+(1-\theta)B)^n(u,v)=\sum_{k=0}^{n}{n \choose k}A^k*B^{n-k}(u,v).$$
\subsection{Proof of Proposition \ref{independent}}
This is a simple consequence of conditional probabilities.
$$P(X+Z<x, Y< y)=\int_{-\infty}^{\infty}P(X+Z<x, Y<y|Z=t)dF_3(z)=$$ $$=\int_{-\infty}^{\infty}P(X+z<x, Y<y)dF_3(z).$$
The second equality uses independence of $(X,Y)$ and $Z$. Using Sklar's theorem,
$$C_5(F_4(x), F_2(y))=\int_{-\infty}^{\infty}H(x-z,y)dF_3(z)=\int_{-\infty}^{\infty}C(F_1(x-z),F_2(y))dF_3(z).$$
Thus, substituting $u=F_4(x)$, $v=F_2(y)$ and $t=F_3(z)$ ends the proof.
\subsection{Details for Example \ref{doublepert}}
The form of the copula when $(X,Y,Z)$ are independent takes into account the range of the variables. In the concrete example of independence (when $C(u,v)=uv$), if the marginal distributions are uniform on $(0,1)$, then $X+Z$ and $Y+Z$ share the same distribution called Irwin-Hall distribution $$F_4(x)=F_5(x)=\begin{cases} 0,& x<0, \\ \frac{1}{2}x^2 & 0\leq x< 1, \\ 1-\frac{1}{2}(x-2)^2 & 1\leq x< 2, \\ 1 & x>2.\end{cases}$$ The joint distribution of $(X+Z, Y+Z)$ is obtained by computing the integral given in Proposition \ref{independent}. The range splits into 8 cases depending on $x$, $y$ and the range of the variable of integration $(0.1)$. Due to symmetry, it is enough to consider 4 cases for $x\leq y$, and flip the result for the other 4 cases by symmetry. The 4 cases are:
\begin{enumerate}
\item $0<x\leq y<1 $, \quad\quad $C_6(F_4(x), F_4(y))=\frac{1}{2}x^2y-\frac{1}{6}x^3$
\item $x< 1 \leq y\leq 1+x$, \quad\quad $C_6(F_4(x), F_4(y))=\frac{1}{2}x^2-\frac{1}{6}(x-y+1)^3$
\item $1+x<y<2$, \quad\quad $C_6(F_4(x), F_4(y))=\frac{1}{2}x^2$
\item $x\leq y, 2>x>1, 2>y>1$, \\
$C_6(F_4(x), F_4(y))=1-\frac{1}{2}(2-x)^2-\frac{1}{6}(2-y)^2(2-y+3x-3).$
\end{enumerate}
Other cases out of this range give value $0$ for the joint probability, $F_4(x)$ or $F_4(y)$, which are irrelevant for the process leading to the copula. Using the inverse of $F_4$ in the formula and substituting $u$ and $v$ in the formula of the range of $x$ and $y$ gives the needed formula of $C_6(u,v)$. Other assertions of this example are obvious.

\subsection{Proof of proposition \ref{bothcomp}}
In general, conditioning on $(Z_1,Z_2)$, we have $$C_7(F_7(x), F_8(y))=P(X+Z_1<x, Y+Z_2<y)=$$ $$=\int_{-\infty}^{\infty}\int_{-\infty}^{\infty}P(X+Z_1<x, Y+Z_2<y|Z_1=z_1, Z_2=z_2)dG_1(z_1)dG_2(z_2).$$ Under the conditions of this proposition, using Sklar's theorem, this formula leads to
$$C_7(F_7(x),F_8(y))=P(X+Z_1<x, Y+Z_2<y)=$$ $$=\int_{-\infty}^{\infty}\int_{-\infty}^{\infty}P(X<x-z_1, Y<y-z_2)dG_1(z_1)dG_2(z_2),$$
$$C(u,v)=\int_{-\infty}^{\infty}\int_{-\infty}^{\infty}C(F_1(F_7^{-1}(u)-z_1), F_2(F_8^{-1}(u)-z_1))dG_1(z_1)dG_2(z_2).$$
Using uniform distributions for $X,Y,Z_1,Z_2$, once more $F_7(x)=F_{8}(x)$. The region of integration splits into 4 blocks in each of which the form of the integrand changes:
\begin{enumerate}
\item $I= \{0\leq x, y\leq 1\}$. For $(x,y)$ in this region, the variables of integration vary from $0$ to $x$ and from $0$ to $y$ respectively because both $x-t_1$ and $y-t_2$ are between $0$ and $1$. 
\item $II= \{0\leq x\leq 1, 1<y\leq 2\}$. Here, the second variable of integration splits into two parts. For one of the parts the integrand is $C(x-t_1, y-t_2)$ on $[y-1, 1]$ and for the other part the integrand turns into $C(x-t_1, 1)$ on $[0, y-1]$.
\item $III= \{0\leq y\leq 1, 1<x\leq 2\}$. Here, the first variable of integration splits into two parts. For one of the parts the integrand is $C(x-t_1, y-t_2)$ on $[x-1, 1]$ and for the other part the integrand turns into $C(1, y-t_2)$ on $[0, x-1]$.
\item $IV= \{0\leq x\leq 1, 1<y\leq 2\}$. Here, both variable of integration split each into two parts, giving rise to 4 different integrals. For one of the parts the integrand is $C(x-t_1, y-t_2)$ on $[x-1, 1]\times[y-1,1],$ the second part has integrand $C(x-t_1, 1)$ on $[x-1,1]\times[0, y-1]$, the third part has integrand $C(1, y-t_2)$ on $[0,x-1]\times[y-1,1]$ and the fourth piece has integrand $C(1,1)=1$ on $[0,x-1]\times[0,y-1]$. This uses the fact that the range of the variables of a copula is $[0,1]$ and $u=F_4(x), v=F_4(y)$ in each of the cases.
\end{enumerate}
Proper computation of the listed integrals leads to the provided formula of $C_7(u,v)$.
\subsection{About Proposition \ref{mcopula}}
Copulas with densities $m_i,$ $i=1,2,3,4$ have continuous partial derivatives at 0 and 1. For any copula with continuous partial derivatives at the given points, the tail dependence coefficients are equal to zero because they are a combination of these derivatives.
$$\lambda_L(m_i)=\lim_{u\to 0}\int_{0}^{u}m_i(u,t)dt+\int_{0}^{u}m_i(t,u)dt=0\quad \mbox{and}$$ 
$$\lambda_U(m_i)=2-\lim_{u\to 1}(\int_{0}^{u}m_i(u,t)dt+\int_{0}^{u}m_i(t,u)dt)=0.$$ 
The above formulas are obtained from the copula representation of the two dependence coefficients, the chain rule and the fact that the limit is equal the derivative in the case of the lower tail dependence and equal to $2$ minus the derivative in the case of the upper tail dependence coefficient. Thus, no copula with continuous derivatives at $0$ exhibits lower tail dependence and no copula with continuous derivatives at $1$ exhibits upper tail dependence.

\subsection{Proof of Theorem \ref{theobeta} and Proposition \ref{Mmixing}}
It has been shown that $$\tilde{C}^n(u,v)=(1-\theta)^nC^{n}(u,v)+(1-(1-\theta)^n)\Pi(u,v).$$ Moreover, by Remark \ref{rbeta}, 
$$\beta(\tilde{C}^n)=(1-\theta)^n\beta_n(C).$$ As $0\leq \beta_n(C)\leq 1$, it follows that under any scenario, 
$\beta_n(\tilde{C})=\beta(\tilde{C}^n)\to 0$ exponentially fast as $n\to\infty$ for any $\theta\in (0,1)$. Therefore the first point of Theorem \ref{theobeta} holds.
The second point of this theorem relies on the fact that when $n\to\infty$, the component containing $M(u,v)$ converges uniformly to $0$.
Thus, $$\lim_{n\to\infty}\beta_n(\hat{C})\leq \lim_{n\to \infty}\sum_{i=1}^{n}{n \choose i}\theta^i(1-\theta)^{n-i}\beta_{i}(C).$$
This inequality uses the triangle inequality to obtain the sum of absolute values before distributing the $\max$ function to obtain the $\beta_i$'s.
This also takes into account the fact that the probability measures induced by the copulas satisfy the relationship that the copulas satisfy.
The conclusion relies on Proposition \ref{sequence}.
The third point of this Theorem is a consequence of Theorem 4 of Longla (2015) and Proposition \ref{Mmixing} is a consequence of Proposition \ref{sequence}.
\subsection{Proof of Proposition \ref{sequence}}
Without loss of generality, we can assume that the limit is $0$. Given that $a_n$ converges to $0$, for any $\varepsilon>0$, assume $N$ is an integer such that for all $n>N$, $|a_n|<\varepsilon/2$.
We have
$$
b_n=\sum_{i=1}^n {n \choose i} \theta^i (1-\theta)^{n-i}a_i=\sum_{i=1}^N {n \choose i} \theta^i (1-\theta)^{n-i}a_i+\sum_{i=N+1}^n {n \choose i} \theta^i (1-\theta)^{n-i}a_i$$ 
$$\mbox{Thus,}\quad |b_n|< \sum_{i=1}^N {n \choose i} \theta^i (1-\theta)^{n-i}|a_i| + \varepsilon/2.$$
This is true because in the second sum, $|a_i|<\varepsilon/2$ for all $i$ and the portion of the binomial sum that remains in the formula is less than 1.
Let $M=\max(|a_i|, 1\leq i\leq N)$. We have
$$|b_n|< M\sum_{i=1}^N {n \choose i} \theta^i (1-\theta)^{n-i} + \varepsilon/2.$$
It is easy to show that $${n \choose i}\leq n^i. $$ 
Thus, $|b_n|< M\sum_{i=1}^N n^i \theta^i (1-\theta)^{n-i} + \varepsilon/2< M(1-\theta)^n\sum_{i=1}^N n^i \theta^i (1-\theta)^{-i} + \varepsilon/2.$ The last sum is a partial sum of a geomeric series, leading to 
$$|b_n|< M(1-\theta)^n (n\theta(1-\theta)^{-1})\frac{1-(n\theta(1-\theta)^{-1})^{N}}{1-(n\theta(1-\theta)^{-1})} + \varepsilon/2.$$
It follows that $$|b_n|<M(1-\theta)^n (-1+(n\theta(1-\theta)^{-1})^{N}) + \varepsilon/2=K (1-\theta)^n n^N+\varepsilon/2,$$ where $K=M(\theta(1-\theta)^{-1})^N$. Thanks to the fact that for any fixed value of $N$ and $0<\theta<1$, $(1-\theta)^n n^N\to 0$ as $n\to \infty$, we can select $N$ such that $\displaystyle K(1-\theta)^n n^N<\varepsilon/2,$ leading to $0<b_n<\varepsilon$ for any $n>N$. Thus, $b_n\to 0$ as $n\to \infty$.

\subsection{About Theorem \ref{phi} and Theorem \ref{psimixing}}
The proof of this result essentially repeats the steps of the proof of Proposition \ref{Mmixing} and Theorem \ref{theobeta}.
The proof of Theorem \ref{psimixing} has two phases. Firstly, it is enough to prove it for a convex combination of two copulas because any convex combination of copulas can be rewritten as a convex combination of two copulas and by mathematical induction the result can be assumed true for any convex combination. Under the conditions of this theorem, if each of the copulas generates $\psi$-mixing, it follows that each of them generates $\phi$-mixing by the properties of mixing coefficients. Thus, by Theorem 4 of Longla (2015), the Markov chains generated by this
convex combination are $\phi$-mixing. $\phi$-mixing implies ergodicity. Thus, it is enough to show here by Theorem \ref{psi} that $\psi_1(a C_1+(1-a)C_2) \ne \infty$. But $\psi_1(a C_1+(1-a)C_2)\leq a\psi_1(C_1)+(1-a)\psi_1(C_2) \ne \infty$ because each of the copulas generates $\psi$
-mixing. This leads to $\psi$-mixing for Markov chains generated by the convex combination of these copulas.  It is good to remark here that when one of the copulas doesn't generate $\psi$-mixing we might have in different cases different results.

\section*{Acknowledgements} {The authors thank the anonymous reviewers and the associate editor who have helped improve the quality of this paper by making valuable comments and suggestions.}

\end{document}